\documentclass {amsart}

\usepackage[T1]{fontenc}
\usepackage {microtype}
\usepackage {amsmath}
\usepackage {mathtools}
\usepackage {thm-restate}
\usepackage {amsfonts}
\usepackage {mathrsfs}
\usepackage {hyperref}
\usepackage {parskip}
\usepackage {dirtytalk}
\usepackage {tikz-cd}
\usepackage {tikz}
\usetikzlibrary {decorations.pathreplacing,angles,quotes}
\usepackage {tgpagella}
\usepackage[citestyle=alphabetic,bibstyle=alphabetic]{biblatex}

\newcommand {\Out}      {\text{Out}}

\newcommand {\CAT}      {\text{CAT}}
\newcommand {\ccat}     {\text{c-cat}}
\newcommand {\pcat}     {\text{p-cat}}
\newcommand {\cat}     {\text{cat}}
\newcommand {\Id}       {\text{Id}}
\newcommand {\asdim}    {\text{asdim}}
\newcommand {\interior}    {\text{int}}
\newcommand {\R}        {\mathbb{R}}
\newcommand {\Z}        {\mathbb{Z}}
\newcommand {\N}        {\mathbb{N}}
\newcommand {\U}        {\mathcal{U}}
\newcommand {\V}        {\mathcal{V}}
\newcommand {\B}        {\mathcal{B}}
\newcommand {\C}        {\mathcal{C}}
\newcommand {\D}        {\mathcal{D}}

\newcommand {\Rplus}    {\mathbb{R}_{+}}
\newcommand {\binfty}   {\mathfrak{B}_{\infty}}
\newcommand {\shr}      {\text {Sh} _{\rho}}
\newcommand {\im}       {\text{Im}}
\newcommand {\Top}       {\textbf{Top}}
\newcommand {\norm} [1] {\left\| #1 \right\|}

\newtheorem   {theorem}               {Theorem}       [section]
\newtheorem   {lemma}       [theorem] {Lemma}
\newtheorem   {claim}       [theorem] {Claim}
\newtheorem   {proposition} [theorem] {Proposition}
\newtheorem   {corollary}   [theorem] {Corollary}
\theoremstyle {definition}
\newtheorem   {definition}  [theorem] {Definition}

\newtheorem   {example}  [theorem] {Example}

\title{On the Coarse Lusternik-Schnirelmann Category of Groups}
\author{Aditya De Saha}
\date{\today}

\address{A. De Saha, Department of Mathematics, University
of Florida, 433 Little Hall, Gainesville, FL 32601, USA}
\email{a.desaha@ufl.edu}

\begin {document}

\begin{abstract}
    We introduce a coarse analog of the classical Lusternik-Schnirelmann category,
    denoted by \( \ccat \), defined for metric spaces in the coarse homotopy
    category. This provides a new tool for studying large-scale
    topological properties of groups and spaces. We establish that \( \ccat \) is
    a coarse homotopy invariant and prove a lower-bound \( \pcat(\Gamma)
    \leq \ccat(\Gamma) \) for geometrically finite groups \( \Gamma \), where \(
    \pcat \) introduced by Ayala and co-authors in 1992. We also prove an upper
    bound \( \ccat(\Gamma) \leq \asdim (\Gamma) \) for bicombable 1-ended groups
    which are semistable at \( \infty \).
\end{abstract}

\maketitle

\section {Introduction}

Numerical invariants proved to be useful in all areas of mathematics, and
perhaps the most popular among them is the concept of dimension. In coarse
geometry dimension appeared due to Gromov in several forms
\cite{gromovAsymptoticInvariantsInfinite1993, 
      gromovPositiveCurvatureMacroscopic1996}.
The most studied among them is the asymptotic dimension
\cite{bellAsymptoticDimension2008} 
which is an important invariant in geometric group theory. Gouliang Yu proved
\cite{yuNovikovConjectureGroups1998}
that the finiteness of the asymptotic dimension of a group $ \Gamma $
implies the Novikov Higher Signature conjecture for $ \Gamma $
and most of the satellite conjectures.

Proving that the asymptotic dimension of a certain group or a class of groups
is finite is often a great challenge. It was proven for hyperbolic
groups 
\cite{gromovHyperbolicGroups1987,roeHyperbolicGroupsHave2005},
nilpotent groups 
\cite{bellAsymptoticDimension2008},
solvable groups 
\cite{dranishnikovAsymptoticDimensionDiscrete2006},
arithmetic groups,
\cite{jiAsymptoticDimensionIntegral2004}
and for mapping class groups
\cite{bestvinaAsymptoticDimensionMapping2010}.
The next challenge are the groups 
$ \Out(F_n) $
and Helly groups 
\cite{bandeltCliqueGraphsHelly1991,chalopinHellyGroups2025}.

We note that Gromov's definition of asymptotic dimension was a translation
of Lebesgue's definition of the covering dimension to the language of coarse
geometry. In this paper we do a similar thing with another numerical invariant
from classical topology, the Lusternik-Schnirelmann category $ \cat(X) $
(LS-category for short). In topology, the LS-category is a lower bound for
dimension, $ \cat(X)\le\dim(X) $, some generalization work has been done in 
\cite{margolisCoarseHomologicalInvariants2024,
srinivasanLusternikSchirelmannCategoryPeano2015}. Our main result is a similar 
inequality
$$
    \ccat(X)\le \asdim(X)
$$
for certain classes of groups and metric spaces
(Theorem \ref{thm:ccat-asdim-grp}, \ref{thm:ccat-asdim}). Thus in the open problems 
about the finiteness of asymptotic dimension perhaps the first step would be to 
try to prove the finiteness of the coarse LS-category.

If a classifying space $ B\Gamma $ of a group $ \Gamma $ is compact, then its
universal cover $ E\Gamma $ with lifted geodesic metric from $ B\Gamma $ is
coarsely equivalent to $ \Gamma $ with a word metric. In that case $ E\Gamma $
is a proper metric space. The LS-category in the category of locally compact
metric spaces and proper maps was defined and studied
\cite{ayalaLusternikSchnirelmannInvariantsProper1992} 
before. We prove the following comparison result in theorem \ref{thm:pcat-ccat}:
\begin{theorem}
    For a geometrically finite group \( \Gamma \) we have the following inequality
    \[
        \pcat(E\Gamma) \leq \ccat(E\Gamma).
    \]
\end{theorem}

\section {The Coarse Homotopy Category}

Let us give a brief description of coarse maps and coarse homotopies, when
restricted to metric spaces. The following can be generalized for general coarse
structures, the interested reader is encouraged to look at
\cite{roeLecturesCoarseGeometry2003, mitchenerCoarseHomotopyGroups2020} for a
more detailed overview.

Throughout the paper we restrict our attention to \emph{proper} metric spaces
(i.e.\ closed balls are compact) which are, in addition, chain-connected in the
following sense.

\begin{definition} \label{def:chain-connected}
    Let \( R > 0 \). A metric space \( X \) is \emph{\( R \)-chain-connected} if
    for every \( x, y \in X \) there is a finite sequence \( x = x_0, x_1, \ldots,
    x_n = y \) in \( X \) with \( d(x_i, x_{i+1}) \leq R \) for all \( i \) (cf.
    \cite[Chapter I.8]{bridsonHaefligerMetricSpacesNonpositive1999}). We call \( X
    \) \emph{chain-connected} if it is \( R \)-chain-connected for some \( R > 0
    \).
\end{definition}

This is a mild assumption, satisfied automatically by geodesic and path-metric
spaces, and by any group with a word metric. Its purpose is to guarantee that our
spaces admit coarse rays: a sequence of \( R \)-jumps escaping every bounded set
is itself a coarse map \( \Rplus \to X \). From now on, all metric spaces are
assumed to be proper and chain-connected, unless stated otherwise.

\begin{definition}
    Let \( X \) and \( Y \) be metric spaces. A (not necessarily continuous)
    function \( f : X \to Y \) is called \emph{controlled}, or \emph{bornologous}
    if for every \( r>0 \) there exists a
    \( S>0 \) such that
    \[
        d(x,x') < r \implies d(f(x), f(x')) < S
    \]
    for all \( x,x' \in X \). The function \( f \) is called \emph{proper} if for
    any bounded subset \( B \subset Y \), the preimage \( f^{-1} (B) \in X \) is
    bounded. \( f \) is called \emph{coarse} if it is both controlled and proper.
\end{definition}

\begin{definition}
    Let \( \rho : \Rplus \to \Rplus \) be any function. A
    function \( f : X \to Y \) between metric spaces is called \emph{\( \rho
        \)-controlled} if it satisfies
    \[
        d( f(x), f(x') ) \leq \rho (d(x,x'))
    \]
    for all \( x,x' \in X \). Clearly, \( f \) is controlled if
    \( f \) is \( \rho \)-controlled for some \( \rho \).
\end{definition}

\begin{definition} \label{def:coarse-embedding}
    Let \( \rho_1, \rho_2 : \Rplus \to \Rplus \) be two functions going to
    infinity, and let \( f : X \to Y \) be a function between metric spaces. \( f
    \) is called a \emph{\( (\rho_1, \rho_2) \)-coarse embedding} if it satisfies
    \[
        \rho_1 (d(x,x')) \leq d( f(x), f(x') ) \leq \rho_2 ( d(x,x'))
    \]
    for all \( x,x' \in X \). 
    A function \( f: X \to Y \) is a \emph{coarse embedding} if \( f \) is a
    \( (\rho_1, \rho_2) \)-coarse embedding for some increasing functions
    \( \rho_1, \rho_2 \) going to infinity.
\end{definition}

\begin{definition}
    Two functions \( f,g : X \to Y \) between metric spaces are \emph{close} or
    \emph{uniformly bounded distance apart} if there exist a constant \( M \in
    \Rplus \) such that \( d(f(x), g(x)) \leq M \) for all \( x \in X \). \( f \)
    is called \emph{bounded} if it is close to a constant map.
\end{definition}

\begin{definition} \label{def:close}
    Two metric spaces \( X,Y \) are called \emph{coarsely equivalent} if there
    exist controlled maps \( f : X \to Y \) and \( g : Y \to X \) such that the
    compositions \( f \circ g \) and \( g \circ f \) are close to \( 1_Y \) and
    \( 1_X \) respectively.
\end{definition}

\begin{example}
    Standard Examples of coarse maps would be
    \begin{itemize}
        \item The floor function \( \lfloor . \rfloor : \R \to \Z \), and the 
              inclusion \( \Z \to \R \). This shows us \( \R \) and \( \Z \) are 
              coarsely equivalent.

        \item In general a similar argument shows us that the integer lattice \(
              \Z^{n} \subset \R^{n} \) is coarsely equivalent to \( \R^{n} \).

        \item Let \( M \) be a complete simply-connected Riemannian manifold of
              non-positive sectional curvature. For a point \( p \in M \), The  
              exponential map \( \exp : T_p M \to M \) is a distance-increasing 
              diffeomorphism. The inverse \( \log : M \to T_p M \) is therefore a 
              coarse map.
    \end{itemize}
\end{example}

\begin{definition} \label{def:coarse-ray}
    Any coarse map \( \Rplus \to X \) is a \emph{coarse ray}.
\end{definition}

\subsection {Coarse Homotopy}

The purpose of this sub-section is to define the notion of homotopy in the
coarse category, as defined in \cite{mitchenerCoarseHomotopyGroups2020}. These
homotopies have to end eventually, but the end time will be allowed to depend
on the given point in the metric space (and to go to infinity as one goes to
infinity). These will be measured by coarse maps \( p : X \to \Rplus \), which
are sometimes called \say{basepoint projections}. Note that if we fix a
basepoint \( x_0 \in X \), there is a natural choice for the projection
\( p_0 : X \to \Rplus \) given by \( p_0 (x) = d(x, x_0) \). \( p_0 \) is also
called the \emph{standard basepoint projection}.

Before defining coarse homotopies, we need to define cylinders.

\begin{definition}
    Let \( X \) be a metric space, and let \( p: X \to \Rplus \) be a coarse map.
    We define the \( p \)-cylinder
    \[
        I_p X = \{ (x,t) \in X \times \Rplus \ | \ t \leq p(x) \}
        .
    \]
  If we have a pointed metric space \( (X,x_0) \), we have a natural choice of
  coarse map, \( p(x) = d(x, x_0) \). In this case we will omit the \( p \) and
  denote \( I(X) = I_p(X) \).
\end{definition}

We have inclusions \( i_0, i_1 : X \to I_p X \) defined by the formulas \( i_0
(x) = (x,0) \) and \( i_1 (x) = (x, p(x) ) \). The canonical projection \( q
: I_p X \to X\) is a coarse map, and identities \( q \circ i_0 = q \circ i_1 =
1_X \) clearly hold.

\begin{definition}
    Let \( X, Y \) be metric spaces. A \emph{coarse homotopy} is a coarse map \( H
    : I_p X \to Y\) for some coarse map \( p : X \to \Rplus \).

    We call coarse maps \( f,g : X \to Y \) \emph{coarsely homotopic} if there is a
    coarse map \( p : X \to \Rplus \) and a coarse homotopy \( H : I_p X \to Y \)
    such that \( H \circ i_0 = f \) and \( H \circ i_1 = g \).

    Let \( f : X \to Y \) be a coarse map between metric spaces. We call the map \(
    f \) a \emph{coarse homotopy equivalence} if there is a coarse map \( g : Y \to
    X \) such that the compositions \( g \circ f \) and \( f \circ g \) are
    coarsely homotopic to the identities \( 1_X \) and \( 1_Y \) respectively.
\end{definition}

We mention a few properties of coarse homotopies, proofs can be found in
\cite{mitchenerCoarseHomotopyGroups2020}.

\begin{itemize}
    \item If two coarse maps \( f,g : X \to Y \) are close (\ref{def:close}), then they 
         are coarsely homotopic.
    \item Coarse homotopy is an equivalence relation.
    \item For path-metric spaces, the choice of the map \( p: X \to \Rplus \) does not 
          matter:
        \begin{lemma} \label{lem:chom-invariant}
              \cite[lemma 2.6]{mitchenerCoarseHomotopyGroups2020} Let \( X \) be a
              path-metric space. For \( x_0 \in X \), let
              \( p_0 : X \to \Rplus \): \( x \mapsto d(x,x_0) \) be the standard
              basepoint projection. Let \( q : X \to \Rplus \) be any coarse map. Then
              any coarse homotopy \( H : I_q X \to Y \) between \( f,g : X \to Y \)
              gives rise to a coarse homotopy \( \bar H : I_{p_0} X \to Y \) between \(
              f \) and \( g \).
          \end{lemma}

\end{itemize}

\begin{definition}
    Two spaces \( X, Y \) are called \emph{coarsely homotopy-equivalent} if there
    exist coarse maps \( f : X \to Y \) and \( g : Y \to X \) such that the
    compositions \( f \circ g \) and \( g \circ f \) are coarsely homotopic to
    \( \Id_Y \) and \( \Id_X \) respectively.
\end{definition}

\begin{example}
    Let \( M \) be a complete simply-connected Riemannian manifold of non-positive
    sectional curvature and let \( p \in M \). The exponential map \( exp : T_p M 
    \to M \) is distance-increasing diffeomorphism, hence the inverse \( log : M
    \to \R^{n} \) is a coarse map. Although \( \exp \) is not coarse (so \( \log
    \) is not necessarily a coarse equivalence), one can construct a \emph{coarse
    homotopy inverse} of \( \log \) by using a radial shrinking function on \(
    \R^{n} \) \cite[Example 2.7]{mitchenerCoarseHomotopyGroups2020}. This proves
    that \( \log \) is a coarse homotopy equivalence. In particular, \( \R^{n} \)
    and the hyperbolic space \( \mathbb{H}^{n} \) are coarsely homotopy
    equivalent.
\end{example}

\begin{definition}
    A metric space \( X \) is \emph{coarsely path-connected} if any two coarse
    rays are coarsely homotopic to each other.
\end{definition}
Note that the notion of \emph{ coarsely path-connectedness} is very similar to the
notion of semistability at \( \infty \) in the group-theory world; we make this
precise in Section \ref{sec:semistability}.

\subsection {Coarse Lusternik-Schnirelmann Category}

One of the classical homotopy invariants in the category \( \Top \) of
topological spaces is the Lusternik-Schnirelmann category \( \cat (X) \), or
LS-category in short. The (reduced) \( \cat (X)\) is defined as the smallest
number \( k \) such that there is an open covering \( \{ U_i \}_{0 \leq i \leq
k} \) of \( X \) with the property that each inclusion map \( U_i
\hookrightarrow X \) is nullhomotopic. We aim to define an analog of this for
the coarse homotopy category.

In the coarse category, instead of being nullhomotopic, we need to define a
notion of sets being \say{small}:

\begin{definition}

    Let \( X \) be a metric space. A subset \( A \subseteq X \) is called
    \emph{coarsely categorical} if there exist coarse maps \( \alpha : \Rplus \to X
    \) and \( j : A \to \Rplus \) such that the following diagram commutes up to
    coarse homotopy:

    \begin{center}

        \begin{tikzcd}
            A \arrow [rd, "j"] \arrow [rr, hook] & & X \\
            & \Rplus \arrow [ur, "\alpha"]
        \end{tikzcd}

    \end{center}

    Where the top horizontal map is the canonical inclusion of \( A \) into \( X
    \).

\end{definition}

Now we can define coarse LS-category:

\begin{definition}[coarse LS-category]
    Let \( X \) be a metric space. The (reduced) coarse LS-category of \( X \),
    denoted by \( \ccat (X) \) is the least number \( k \) such that there exists
    a covering \( \{U_i\} _{0 \leq i \leq k} \) of \( X \) by \( k+1 \) coarsely
    categorical sets. If no finite covering of \( X \) by coarsely categorical
    sets exists then \( \ccat (X) = \infty \).
\end{definition}

Note that if \( X \) admits no coarse ray at all (for instance if \( X \) is
bounded), then no non-empty subset of \( X \) can be coarsely categorical, since
the map \( \alpha \) in the definition above would itself be a coarse ray.
Hence \( \ccat(X) = \infty \) for such \( X \), by the convention above. Our
standing assumption that spaces be proper and chain-connected
(Definition \ref{def:chain-connected}) rules this out for every unbounded space
we consider.

In \cite{mitchenerCoarseHomotopyGroups2020} the authors proved that \( \pi_0
(\Rplus) \) is trivial, so any coarse map \( \Rplus \to \Rplus \) is coarsely
homotopic to the identity. As a result, we can say that for a metric space \( X
\), \( \ccat(X) = 0 \) if and only if \( X \) is coarsely homotopy equivalent to
\( \Rplus \). 

The following two results show that \( \ccat \) is a coarse homotopy invariant.

\begin{lemma}
  \label{lem:ccat-dom}
    If there exist coarse maps \( f : X \to Y \) and \( g : Y \to X \) such that
    \( f \circ g : Y \to Y \) is coarsely homotopic to the identity \(
    1_Y \), then \( \ccat (Y) \leq \ccat (X) \).
\end{lemma}

\begin{proof}
    Let \( \ccat X = k \), and 
    \( \{ U_i \} _{i=0} ^{k} \) be a covering of \( X \) by coarsely categorical
    sets. Let \( V_i = g^{-1} (U_i) \) be their pre-images in \( X \).
    Suppose \( G : I_p(Y) \to Y \) is a coarse homotopy between \( \Id_Y \) and
    \( f \circ g \), and \( H_i : I_q (U_i) \to X \) be a coarse contraction of
    \( U_i \) to a ray \( \alpha : \Rplus \to X \).

    The restriction of \( G \) to \( I_p (V_i) \) is a coarse homotopy of the
    inclusion \( V_i \to Y \) to \( f \circ g \). But since \( g(V_i) \subseteq
    U_i \), we can take the composition \( f \circ H \) which is a homotopy
    contracting \( f|_{U_i} \) to the ray \( f \circ \alpha \). Hence composing
    the two homotopies \( G|_{I_p(V_i)} \) and \( f \circ H \), we get a
    homotopy contracting \( V_i \).

    Since \( V_i \)'s are coarsely contractible, \( \ccat Y \leq k \).
\end{proof}

\begin{proposition}
    If \( X, Y \) are metric spaces which are coarsely homotopy equivalent, then
    \( \ccat (X) = \ccat (Y) \).
\end{proposition}

\begin{proof}
  This follows from the above lemma \ref{lem:ccat-dom}. If \( X, Y\) are
  coarsely equivalent, then there are coarse maps \( f: X \to Y \) and \( g : Y
  \to X\) such that \( f \circ g \) and \( g \circ f \) are coarsely homotopic
  to identity maps. Using lemma \ref{lem:ccat-dom} we can say \( \ccat (X) \leq
  \ccat (Y)\) and \( \ccat(Y) \leq \ccat (X) \), therefore both quantities are
  equal.
\end{proof}

Since \( \ccat \) is a coarse homotopy-invariant, we can define \(
\ccat(\Gamma) = \ccat(X) \) for any proper geodesic metric space \( X \) on which 
\( \Gamma \) acts properly and cocompactly via isometries.

\begin{example}
    \begin{enumerate}
        \item \( \ccat(X) = 0 \) if and only if \( X \) is coarsely homotopy
              equivalent to \( \Rplus \).
        \item \( \ccat(\R^{n}) = 1 \) for all \( n > 0 \). We can cover \( \R^{n} \)
              with two halves \( A= [0, \infty) \times \R^{n-1} \) and \( B = [0,
              \infty) \times \R^{n-1} \). \( A \) can be coarsely deformed to the
              positive \( x \)-axis, and \( B \) the negative \( x \)-axis. 
              \( \ccat( \R^{n} ) > 0 \) because of proposition \ref{prop:pcat-ccat-ineq}.
        \item Because of the above, \( \ccat (M) = 1 \) for any complete
              simply-connected Riemannian manifold of non-positive sectional 
              curvature, and \( \ccat(\Z^{n}) = 1 \).
        \item For the infinite binary tree \( T_2 \), no two distinct geodesic rays
              are coarsely homotopic. Hence \( \ccat(T_2) = \infty \).
    \end{enumerate}
\end{example}

\section {Coarse LS-category vs proper LS-category}

In \cite{ayalaLusternikSchnirelmannInvariantsProper1992} the authors explored
the idea of proper LS category \( \pcat \), which is quite similar to our notion
of \( \ccat \), where instead of requiring maps to be coarse, they require all
maps to be proper and continuous. This is slightly less restrictive, so as we'll
see from the discussion below, for a lot of reasonable spaces \( X \), we have the
inequality \( \pcat(X) \leq \ccat(X) \).

Before we go into the definitions, let us first prove a lemma that allows us to
\say{upgrade} coarse (not necessarily continuous) maps to proper continuous
maps.

\begin{lemma}
    \label{lem:coarse-to-proper}
    Suppose we have a coarse map \( f : X \to Y \) between metric spaces \( X,Y \)
    with the properties that \( X \) is a finite dimensional simplicial complex
    whose simplices are uniformly bounded, and \( Y \) is uniformly contractible.
    Then there exists a coarse continuous map \( g : X \to Y \) which is 
    close to \( f \).
\end{lemma}

Recall that a metric space \( X \) is \emph{uniformly contractible} if for every
\( R > 0 \) there exists an \( S > R \) such that every \( R \)-ball in \( X \)
is contractible inside a bigger \( S \)-ball in \( X \).

\begin{proof}
    First consider \( X^{(0)} \), the zero-skeleton of \( X \). We can define \(
    g^{0} \) on \( X^{(0)} \) to be the restriction of \( f \), and it is
    continuous and coarse (since \( X^{(0)} \) is discrete). Inductively, suppose
    we have constructed a continuous function \( g^{k}: X^{(k)} \to Y \) on the \(
    k \)-skeleton of \( X \) which is coarse and close to \( f \).

    Using the uniform contractibility of \( Y \), we can extend the map \( g \) to
    a larger map \( g^{k+1} \) as shown by the dashed arrow:

    \begin{equation}
        \begin{tikzcd}
            X^{(k)} \arrow[r, "g^{k}"] \arrow[d] & Y \\
            X^{(k + 1)} \arrow[ur,dashed, "g^{k+1}" below] 
        \end{tikzcd}
    \end{equation}

    Locally the construction is done as follows. Suppose we have a simplex
    \( \varphi : \Delta ^k \to X \) in \( X \). Since \( f \) is coarse, and in
    particular controlled, the images of all simplices under \( f \) are also of
    some bounded diameter. By uniform contractibility of \( Y \),
    the image of \( f \circ \varphi \) is contractible inside some \( K
    \)-neighborhood in \( Y \). The contraction can also be seen as a map of the cone 
    \( \Phi : \Delta ^ {k+1} = C(\Delta ^k) \to Y \)
    which agrees with \( \varphi \) on the \( k \)-skeleton \( \Delta^k \).
    Furthermore, \( \im \Phi \) is contained in \( K \)-neighborhood of \( \im
    \varphi \).

    Since these maps agree on the \( k \)-skeleton, we can paste these maps to
    construct the \( k+1 \)-skeleton \( g^{k+1} \) which is close to \( g^k \).

    Because \( X \) is a finite dimensional simplicial complex, this process
    terminates eventually, and we get our continuous map \( g = g^{n} \) where
    \( n= \dim X \).
\end{proof}

Now we can introduce the notion of proper LS-category as done in
\cite{ayalaLusternikSchnirelmannInvariantsProper1992}. We will be working in
the category \( \binfty \) of non-compact \( T_2 \)-locally compact spaces, and
proper continuous maps as morphisms. Let \( X \in \binfty \) be such a space.

Define a closed subset \( C \subseteq X \) to be \emph{properly deformable to
    \( \Rplus\)} if there exists a diagram in \( \binfty \):
\begin{equation*}
    \begin{tikzcd}
        C \arrow[r, hook] \arrow[d, "r"] & X \\
        \Rplus  \arrow[ur, bend right, "\alpha" right]
    \end{tikzcd}
\end{equation*}
which commutes up to proper homotopy.

\begin{definition}
    \label{def:p-cat}
    Given a space \( X \) in \( \binfty \), \( A \subseteq X \) is said to be
    properly categorical in \( X \) if there is a closed neighborhood of \( A \)
    properly deformable to \( \Rplus \) in \( X \).

    An open covering \( \{ U_{\alpha} \} \) of \( X \) is said to be properly
    categorical if each \( U_{\alpha} \) is properly categorical in \( X \).

    The proper Lusternik-Schnirelmann category (or proper LS category) of \( X
    \), denoted by \( \pcat(X) \) is the least number \( n \) such that \( X \)
    admits a properly categorical open covering with \( n+1 \) elements. If no
    finite properly categorical covering exists then \( \pcat(X) = \infty \).
\end{definition}

This is the \emph{reduced} version of the proper LS-category defined by Ayala et
al.\ in \cite{ayalaLusternikSchnirelmannInvariantsProper1992}, shifted by one so
that it agrees with the reduced convention we use for \( \ccat \). In particular
Proposition 1.6(ii) of \cite{ayalaLusternikSchnirelmannInvariantsProper1992}
translates to: \( \pcat(X) = 0 \) if and only if \( X \) is properly homotopic to
\( \Rplus \).

Using our lemma \ref{lem:coarse-to-proper} we can compare our notions of \(
\ccat(X) \) and \( \pcat(X) \) when our space \( X \) is nice:

\begin{proposition} \label{prop:pcat-ccat-ineq}
    Let \( X \) be a uniformly contractible finite dimensional simplicial complex,
    with a metric such that the simplices are of bounded diameter. Then
    \[
        \pcat(X) \leq \ccat(X).
    \]
\end{proposition}

\begin{proof}
    Suppose \( \ccat(X) = n \), and let \( \{ U_1, U_2, \cdots, U_{n+1} \} \) be a
    coarsely categorical cover of \( X \). Define \( V_i \) to be the smallest 
    simplicial subcomplex of \( U_i \) in \( X \). Since simplices in \( X \) are of
    bounded diameter, and \( X \) is uniformly contractible, we can extend the
    coarse homotopies of \( U_i \) to \( V_i \)'s, and therefore \( \{ V_1, \cdots,
    V_{n+1} \} \) is a coarsely categorical cover of \( X \) as well.

    We can give a simplicial complex structure to each \( I(V_i) \) such that the
    simplices are uniformly bounded. Consider the coarse homotopy \( H_i : I(V_i)
    \to X \). Using lemma \ref{lem:coarse-to-proper} we can say that each \( V_i
    \) is properly categorical. Hence \( \{ \interior(V_1) , \cdots , \interior(V_{n+1}) \} \) gives us
    a properly categorical open cover of \( X \), so \( \pcat(X) \leq n \).
\end{proof}

A good class of examples of such simplicial complexes are universal covers of
groups which have a finite \( K(G,1) \)-complex. These are the so-called
\emph{geometrically finite} groups. Hence we have:

\begin{theorem}
    \label{thm:pcat-ccat}
    For a geometrically finite group \( \Gamma \), we have the inequality
    \[
        \pcat (\Gamma) = \pcat(E \Gamma) \leq \ccat (E \Gamma) = \ccat(\Gamma)
        .
    \]
\end{theorem}

\begin{example}
    \begin{enumerate}
        \item \( \pcat (\R^{n}) = 1 \) for all \( n > 0 \).
        \item Let \( X \) be the euclidean plane without a strip around the negative
            \( x \)-axis: \( X = \R^{2} \setminus (-\infty, 0) \times (-1,1) \). 
            Then \( X \) is properly deformable to the positive \( x \)-axis, so \(
            \pcat (X) = 0 \). But \( X \) is coarsely equivalent to the plane, hence
            \( \ccat (X) = 1 > \pcat (X) \).
        \item Let \( T \) be an embedding of the binary tree \( T_2 \) inside \( [0,1] 
            \times \Rplus \). \( \pcat (T) = \infty \), but since \( T \) is coarsely
            equivalent to \( \Rplus \), \( \ccat(T) = 0 < \pcat(T) \). \( T \)
            with the subspace metric here is not uniformly contractible, so
            proposition \ref{prop:pcat-ccat-ineq} doesn't hold.
    \end{enumerate}
\end{example}

\section {Combable Spaces}

In this section we will define the notion of combable spaces. Most of our
results in the subsequent sections will be about these spaces.
Interested readers can check \cite{katoAsymptoticLipschitzMaps2000} for a
more detailed description of combability and bicombability.

Let \( X \) be a metric space. For convenience, from now on \( \N \) will
represent the set of nonnegative integers.

\begin{definition}
    (\cite{engelCoronasProperlyCombable2021}, definition 2.4)
    By a \emph{combing} on \( X \) starting at  a point \( p \in X \) we mean a
    map
    \[
        C : X \times \N \to X
    \]
    such that
    \begin{enumerate}
        \item \( C(x,0) = p = C (p,n) \) for all \( x \in X \) and \( n \in \N \).
        \item for each bounded subset \( K \subset X \) there is an \( N \in \N \) such that
              for all \( n \geq N \) and \( x \in K \) we have \( C(x,n) = x \).
        \item \( C \) is a controlled map.
    \end{enumerate}
    Moreover, \( C \) is called \emph{proper} if for any bounded subset \( K
    \subset X \) there is a bounded subset \( L \subset X \) and an \( N \in \N \)
    such that \( C ^{-1} (K) \subset L \) for all \( n \geq N \).
\end{definition}

If a given space can be equipped with a combing then it is called
\emph{combable}.

In a similar vein, one can define a bicombing:

\begin{definition}
    Let \( X \) be a metric space. A \emph{bicombing} on \( X \) is a map
    \[
        C : X \times X \times \N \to X ;
        \ \ \ C_p(x, n) \vcentcolon=  C(p, x, n)
        \ \ \ \forall p,x \in X, n \in \N
    \]
    which satisfies
    \begin{enumerate}
        \item \( C_p(x, 0) = C_p(p, n) = p \) for all \( p,x \in X \), \( n \in \N
              \).
        \item For each \( p \in X \) and \( K \subset X \) bounded, there exist \( N \in \N
              \) such that for all \( n \geq N \) and \( x \in X \), \( C_p(x,n) = x \).
        \item \( C \) is a controlled map.
    \end{enumerate}
\end{definition}

One can regard \( C(x, - ) \) or \( C_p(x, - ) \) as a path from the point \( p
\) at time \( 0 \) to the point \( x \) which eventually becomes constant
(encapsulated in point 2). The map \( C \) being \emph{controlled} implies
that these paths are uniformly coarse, and they have the so-called
\say{fellow-travelling property}. If \( G \) is a finitely generated group endowed
with the word metric, then this notion of a (bi)combing is
equivalent to the usual notion used in geometric group theory literature (they
are sometimes called synchronous (bi)combings or also bounded (bi)combings).

A large class of examples of bicombable spaces are \( \CAT(0) \) spaces
\cite{bridsonHaefligerMetricSpacesNonpositive1999}, via the straight-line
combing towards a fixed basepoint. For groups, every biautomatic group is
bicombable \cite{epsteinWordProcessingGroups1992}; this class includes
hyperbolic groups \cite{gromovHyperbolicGroups1987}, virtually abelian groups,
Artin groups of finite type \cite{charneyArtinGroupsFiniteType1992}, and
systolic groups \cite{januszkiewiczSimplicialNonpositiveCurvature2006}. More
generally, every semihyperbolic group
\cite{alonsoBridsonSemihyperbolicGroups1995} is bicombable, a class which also
contains all \( \CAT(0) \) groups and, by a result of Durham, Minsky and Sisto,
mapping class groups \cite{durhamMinskySistoStableCubulations2023}. Bicombable
groups are finitely presented, so a lot of examples of groups that are not
bicombable are groups that are not finitely presentable.

If we consider the singleton bounded set \( K = \{ x \} \) for property 2, there
exist a smallest \( N =\vcentcolon N_p(x) \) such that \( C_p(x,n) = x \) for
all \( n \geq N \). So each path \( C_p(x, -) \) is of \say{length} \( N_p(x) \).

We prove a \say{coarse} version of uniform contractibility for geodesic metric
spaces which are combable.

\begin{definition} [\( \ast \)]
    \label{def:property-p}
    A metric space \( X \) has property \( (\ast) \) if there exist \( \rho :
    \Rplus \to \Rplus \) going to infinity, such that for any \( r \)-ball \( B
        _{p} (r) \) around any point \( p \in X \), we have a (not necessarily
    continuous) map
    \[
        h: I(B _{p} (r)) \to X
    \]
    such that
    \begin{enumerate}
        \item \( h \circ i_0 (x) = x \) and \( h \circ i_1 (x) = p \) for all \(
            x \in B _{p} (r) \).
        \item \( h \) is \( \rho \)-controlled.
    \end{enumerate}
\end{definition}

One can interpret definition \ref{def:property-p} as a coarse version of uniform
contractibility, where the balls are uniformly coarsely contractible.

\begin{lemma}
    \label{lem:bicomb-property-p}
Let \( X \) be a bicombable path-metric space. Then \( X \) has \( (\ast) \).
\end{lemma}

\begin{proof}
  Suppose our space \( X \) has a bicombing
  \[ 
    C : X \times X \times \N \to X 
  .\]
  Because of the reasons discussed above, the map \( \gamma : X \to \Rplus \)
  defined by \( \gamma(x) = N_p(x) \) is a coarse map. We can now define a
  coarse homotopy
  \[
    h' : I_{\gamma} (B_p(r)) \to X
  \]
  \[
    h'(x, t) = C_p (x, \lfloor N_p(x) - t\rfloor)
  \]
  where \( \lfloor . \rfloor \) denotes the floor function. Our function \(
  h'\) has the property \( h' \circ i_0 (x) = C_p(x, N_p(x)) = x \), and \( h'
  \circ i_1 (x) = C_p(x, 0) = p \). Thus property 1 is satisfied. Property 2
  follows from the fact that \( C \) is controlled.

  Since \( X \) is a path-metric space, by lemma \ref{lem:chom-invariant} we get
  the result.
\end{proof}

For the rest of this section, we will consider cases where proper maps can be
\say{upgraded} to coarse maps, so we get another upper bound for \( \ccat \).
First, we need to define the so-called \emph{shrinking map}.

\begin{definition}
  Suppose we have a pointed metric space \( (X,p) \) which admits a combing \( C :
  X \times \N \to X \) for the basepoint \( p \). Suppose \( \rho : \N \to
  \N \) is any map such that \( \rho(t)  \leq t \) for all \( t \in \N
  \). Then one can define the \emph{shrinking map} \( \shr : X \to X \) as
  \begin{equation}
    \label {eqn:shrinking}
    \shr (x) = C \left( x, \rho (N_x) \right)
  \end{equation}
\end{definition}

Here (and from here onwards) will use the notation \( \norm x = d (x, p) \) for
any point \( x \in X \), whenever the base point \( p \) is obvious from
context. 

Note that in such a case, for each \( x \in X \) one can go from \( x \) to \(
\shr (x) \) \say{along the combing \( C \)}, so we get the following lemma:

\begin{lemma} \label{lem:shrinking}
    For a combable pointed metric space \( (X, p) \) and a coarse map \( \rho : \N
    \to \N \) such that \( \rho (t) \leq t \) for all \( t \in \N \), the
    shrinking map \( \shr \) defined in equation \ref{eqn:shrinking} is coarsely homotopic
    to the identity map.
\end{lemma}

% \adi{Properly writing down this proof will be kind of annoying, but is doable.
% Should I write it?}

Using this lemma, we can prove the following generalization of an earlier result
by Thomas Weighill \cite{weighillLiftingCoarseHomotopies2019}:

\begin{proposition} \label{prop:proper-upgrade}
  Let \( (X, p) \) be a combable pointed metric space, and \( Y \) be any proper
  metric space. If two continuous coarse maps \( f,g : X \to Y \) are properly
  homotopic, then \( f,g \) are coarsely homotopic.
\end{proposition}

\begin{proof}
  Let \( h : X \times [0,1] \to Y \) be a proper homotopy from \( f \) to \( g
  \). We can extend \( h \) by defining a new function \( h' : I_p X \to Y \):

  \begin{equation}
    h' (x,t) = \begin {cases}
      h \left(
          x, \frac {t} {\norm x} 
        \right), & x \neq p \\
        p & x = p
    \end {cases}
  \end{equation}

  We will now construct a decreasing function \( \rho : \Rplus \to \Rplus \).
  Since \( Y \) is a proper metric space and \( h' : I_pX \to Y \) is a
  continuous map, we can find integers \( L_k \) such that
  \[
       d ((x,t), (x', t')) \leq \frac {1} {L_k}
       \implies
       d (h'(x,t), h'(x',t')) \leq 1
  \]
  for \( x \in B _{p} (k) \) and \( t \in [0,k] \). Also we can choose \( L_k
  \)'s to be increasing, and all greater than one.

Define \( \rho: \N \to \N \) as the unique map that maps \( \{0, \cdots, L_1 \}
\) \say{linearly} to \( \{0,1\} \), \( \{L_1, \cdots, L_1 + 2L_2\} \) to \(
\{1,2\} \), \( \{L_1 + 2L_2, \cdots, L_1 + 2L_2 + 3L_3 \} \) to \( \{2,3\} \)
    and so on. Figure \ref{fig:rho} describes the map.

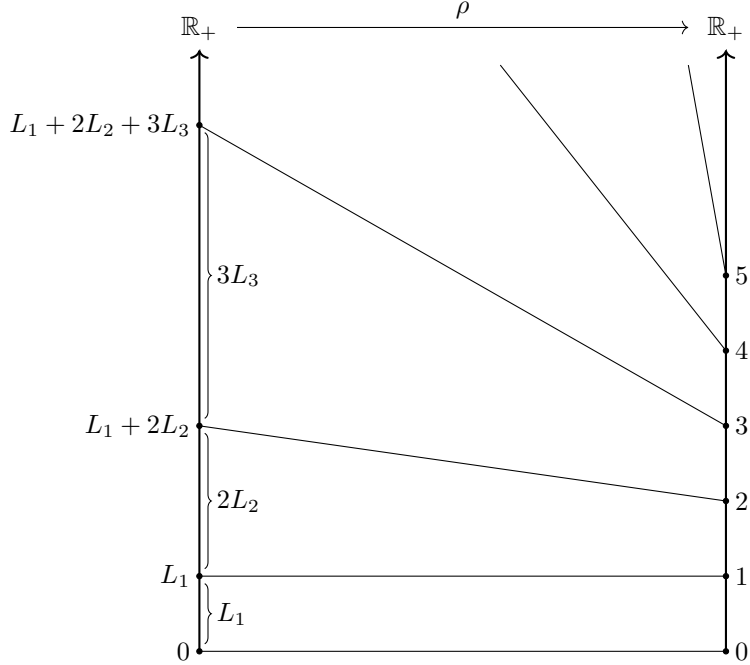
\begin{figure}
    \centering
    \begin{tikzpicture} %% [[[

        %% axes
        \draw [thick, ->] (0, 0) -- (0,8);
        \draw [thick, ->] (7, 0) -- (7,8);

        %% points on left axis
        \filldraw (0,0) circle[radius=1pt] node[anchor=east]{\( 0 \)};
        \filldraw (0,1) circle[radius=1pt] node[anchor=east]{\( L_1 \)};
        \filldraw (0,3) circle[radius=1pt] node[anchor=east]{\( L_1 + 2L_2 \)};
        \filldraw (0,7) circle[radius=1pt] node[anchor=east]{\( L_1 + 2L_2 +
                3L_3 \)};

        %% differences on the left side
        \draw [decoration={brace,mirror,raise=2pt}, decorate]
        (0,0.1) -- (0,0.9) node[midway, right=3pt] {\( L_1 \)};
        \draw [decoration={brace,mirror,raise=2pt}, decorate]
        (0,1.1) -- (0,2.9) node[midway, right=3pt] {\( 2L_2 \)};
        \draw [decoration={brace,mirror,raise=2pt}, decorate]
        (0,3.1) -- (0,6.9) node[midway, right=3pt] {\( 3L_3 \)};

        %% points on right axis
        \filldraw (7,0) circle[radius=1pt] node[anchor=west]{ \( 0 \)};
        \filldraw (7,1) circle[radius=1pt] node[anchor=west]{ \( 1 \)};
        \filldraw (7,2) circle[radius=1pt] node[anchor=west]{ \( 2 \)};
        \filldraw (7,3) circle[radius=1pt] node[anchor=west]{ \( 3 \)};
        \filldraw (7,4) circle[radius=1pt] node[anchor=west]{ \( 4 \)};
        \filldraw (7,5) circle[radius=1pt] node[anchor=west]{ \( 5 \)};

        %% lines between axes
        \draw (0,0) --  (7,0);
        \draw (0,1) --  (7,1);
        \draw (0,3) --  (7,2);
        \draw (0,7) --  (7,3);
        \draw (4,7.8) --  (7,4);
        \draw (6.5,7.8) --  (7,5);

        %% labels for the axes
        \draw (0, 8) -- (0, 8) node[anchor=south] {\( \Rplus \)};
        \draw (7, 8) -- (7, 8)  node[anchor=south] {\( \Rplus \)};
        \draw [->] (0.5, 8.3) -- (6.5,8.3)
        node [midway, anchor=south] {\( \rho \)};

    \end{tikzpicture} %% ]]]
    \caption {A visual description of the map \( \rho : \N \to \N \)}
    \label{fig:rho}
\end{figure}

Having defined \( \rho \), let us define \( H : I_p X \to Y \) as
\[
  H(x,t) = h' \left(\shr(x), t \frac {\rho (\norm x)} {\norm x}\right)
\]

Now we can check that
\begin{claim}
  \label{claim:technical}
  \[
    d( (x,t), (x', t')) \leq 1
    \implies
    d(H(x,t), H(x', t')) \leq 1
  \]
  for all \( (x,t), (x',t') \in I_p X \).
\end{claim}
Assuming the claim, it is easy to see that \( H \) is a coarse homotopy between
\( H \circ i_0 \) and \( H \circ i_1 \). But
\[
  H \circ i_0 (x) = h' (\shr (x), 0) = f(\shr (x)) = (f \circ \shr) (x)
\]
and
\begin{equation*}
  \begin {aligned}
  H \circ i_1 (x) & = h' (\shr (x), \rho (\norm x)) \\
                  & = h' (\shr (x), \norm {\shr(x)}) \\
                  & = g (\shr (x)) = (g \circ \shr ) (x).
  \end {aligned}
\end{equation*}
So our constructed \( H \) is a coarse homotopy between \( f \circ \shr \) and
\( g \circ \shr \). Using lemma \ref{lem:shrinking}, we get our required result.

It remains to prove Claim \ref{claim:technical}. Let's assume we are using the
supremum metric on \( I_p X \); the proof for other metrics is similar.

  Suppose \( d((x,t), (x',t')) \leq 1 \). Therefore \( d(x,x') \leq 1 \) and \(
  d(t,t') \leq 1\). We can assume both points \( (x,t) \) and \( (x',t') \) lie
  in
  \begin{equation*}
    \{ (x,t) \in I_p X \ | \ 
    L_1 + \cdots + (k-1) L_{k-1} \leq \norm x \leq
    L_1 + \cdots + k L_k\}
  \end{equation*}
  for some \( k \in \N \). For convenience, define \( y = \shr(x), y'= \shr(x'),
  s = t \frac {\rho (\norm x)} {\norm x}, s' = t' \frac {\rho (\norm x')} {\norm
    x'}\). Therefore, \( d(x, x') \leq 1 \) implies
  \begin{equation} \label{eqn:inequality1}
    d (y, y') = d (\shr (x), \shr(x')) \leq \frac {1} {k L_k}
       \leq \frac {1} {L_k}.
  \end{equation}
  (Even if both points \( (x,t), (x',t') \) do not lie in such a set, we can take
  the smallest \( k \) such that one of the points lie in such a set, and the
  inequality is still true.)
  Similarly one can argue that (assuming \( \norm x \leq \norm {x'} \))

  \begin{equation} \label{eqn:inequality2}
    \begin {aligned}
    d (s, s') = 
      d \left(
        t \frac {\rho (\norm x)} {\norm x},
        t' \frac {\rho (\norm x')} {\norm x'},
      \right)
        & \leq \frac {\rho (\norm x)} {\norm x} d(t, t') \\
        & \leq \frac {k} {L_1 + \cdots + k L_k} d (t, t') \\
        & \leq \frac {1} {L_k} d(t, t') \leq \frac {1} {L_k}.
    \end {aligned}
  \end{equation}

  Thus if \( x, x' \) satisfy \( L_1 + \cdots + (k-1) L_{k-1} \leq \norm {x}
  \leq L_1 + \cdots + k L _k \), then combining inequalities
  \ref{eqn:inequality1} and \ref{eqn:inequality2} we get 
  \begin{equation*}
    \begin {aligned}
     d((x,t), (x', t')) \leq 1 
       & \implies d((y, s), (y', s')) \leq \frac {1} {L_k} \\
       & \implies d(h'(y, s), h'(y', s')) \leq 1 \\
       & \implies d(H(x, t), H(x', t')) \leq 1 \\
    \end {aligned}
  \end{equation*}
  This is independent of \( k \), so is true for all choices of \( x,x', t,
  t'\). Hence proved.
\end{proof}

Now using our proposition \ref{prop:proper-upgrade} we can now establish the
following upper bound for \( \ccat \):

\begin{theorem}
  \label{thm:ccat-pcat}
  Let \( (X,p) \) be a pointed metric space. Then \( \ccat (X) \leq k\), where
  \( k \) is the least number such that \( X \) can be covered by \( k+1 \)
  subsets each of which admits a combing, and is properly contractible to a
  proper ray.
\end{theorem}

\section {Semistability at Infinity} \label{sec:semistability}

Firstly recall the notion of a Freudenthal end of a space \( X \). The
\emph{space of ends of X}, denoted by \( \mathscr{F}(X) \) is the inverse limit
\( \mathscr{F} (X) = \varprojlim \pi_0 (X-K) \) where \( K \) ranges over the family of
compact subsets of \( X \), and \( \pi_0 \) stands for the set of connected
components. Roughly speaking, these are the ``Connected components at
infinity''.
If our space is a group \( \Gamma \) with the word metric, it is known that \(
\mathscr F (\Gamma)\)  can have either \( 0, 1, 2 \) or \( \infty \) elements.

We recall the following definition of semistability at infinity, used here
\cite{mihalikSemistabilityInfinitySimpleConnectivity1996}.

\begin{definition}
    A set \( S \) in a space \( X \) is \emph{unbounded} if \( S \) is contained
    in no compact subset of \( X \). Two proper rays \( r,s : \Rplus \to X \)
    \emph{converge to the same end} of \( X \) if for any compact set \( C \) in
    \( X \), there exists a number \( N(C) \) such that \( r([N,\infty)) \) and
    \( s([N,\infty)) \) are subsets of the same unbounded component of \( X - C
    \). A connected space \( X \) is called \emph{semistable at \( \infty \)} if
    any two proper rays that converge to the same end of \( X \) are properly
    homotopic, i.e.\ there exists a proper map \( H : \Rplus \times [0,1] \to X
    \) such that \( H(t,0) = r(t) \) and \( H(t,1) = s(t) \).
\end{definition}

A finitely generated group \( \Gamma \) is \emph{semistable at \( \infty \)} if
its Cayley graph, equipped with the word metric associated to some (equivalently
any) finite generating set, is semistable at \( \infty \) as a space. It is known
that any \( 0 \) or \( 2 \)-ended group is semistable at \( \infty \), but the
other two cases are open problems.

\begin{claim} \label{clm:groups-semistable-pi0}
For a finitely presented discrete group \( \Gamma \) with the word metric, \( \Gamma \)
is 1-ended semistable at \( \infty \) if and only if \( \Gamma \) is coarsely
path-connected.
\end{claim}

\begin{proof}
    Consider \( X \) to be the Cayley complex of the group \( \Gamma \) with
    some generating set, which is coarsely equivalent to \( \Gamma \).

    Suppose \( \Gamma \) is \( 1 \)-ended and semistable at \( \infty \). Let \(
    \alpha, \beta : \Rplus \to X \) be two coarse rays. By lemma
    \ref{lem:coarse-to-proper} we can take \( \alpha, \beta \) to be coarse
    continuous maps. Since \( \Gamma \) is semistable at \( \infty \), \( \alpha,
    \beta \) are properly homotopic. By proposition \ref{prop:proper-upgrade},
    \( \alpha, \beta \) are coarsely homotopic.

    For the converse, suppose \( \Gamma \) is coarsely path-connected. Clearly \(
    \Gamma \) is \( 1 \)-ended. Let \( p,q : \Rplus \to X \) be proper rays. Up to
    reparametrization, we can take \( p,q \) to be coarse continuous. Since \( X
    \) is coarsely path-connected, there exists a coarse homotopy \( H_c :
    I(\Rplus) \to X \) between \( p \) and \( q \). We can give \( I(\Rplus) \) a
    uniformly contractible simplicial complex structure, so we can use lemma
    \ref{lem:coarse-to-proper} to get a continuous coarse map \( H' : I_p(\Rplus)
    \to X \). Using \( H' \), we get a proper homotopy between \( p \) and \( q
    \).
\end{proof}

\section {Asymptotic dimension and Coarse LS-Category}

The main result of this section is the following:
\begin{theorem}
    For a bicombable proper geodesic metric space \( X \) which is coarsely
    path-connected, we have the inequality \(
    \ccat(X) \leq \asdim (X) \).
\end{theorem}

For the majority of this section and what follows it, we will be considering
based metric spaces \( (X, p) \), i.e., a metric space \( X \) and a basepoint
\( p \in X \). But before we can prove the theorem, we need to develop some
machinery.

\subsection {Dispersed sets and families} 

Let's define some notions of \say{dispersed} sets, defined in
\cite{bellAsymptoticDimensionBedlewo2005}.

From here on, we will use the notation
\[
    B_p(r) = \{ x \in X \ | \ d(x,p) < r \}
\]
and
\[
    A_p(r,R) = B_p(R) \setminus B_p(r)
\]
for any metric space \( X \), \( p \in X \) and real numbers \( R > r > 0 \).

\begin{definition}
    \label {def:dispersed-sets}
    Given a based metric space \( (X, p) \) and a discrete subset \( U \subset X
    \), \( U \) is called \emph{dispersed} if the function
    \[
        \partial (r) = \inf \{ d(x_1, x_2) : x_1 \neq x_2 \in U \setminus B _{p} (r) \}
    \]
  goes to infinity (\( lim _{r \to \infty} \partial (r) = \infty \)).
\end{definition}

\begin {definition}
    \label{def:dispersed-family}
    Given a based metric space \( (X,p) \) and a family of subsets \( \U \) of
    \( X \), we say the family is dispersed if
    \begin{enumerate}
      \item Elements of \( \U \) are bounded and pairwise disjoint.
      \item The function
        \[
          \partial (r) = \inf \{ d(D_1, D_2)
          \ | \ d(p, D_i) > r \text { and } D_1 \neq D_2 \in \U \}
        \]
        goes to infinity (\( lim _{r \to \infty} \partial (r) = \infty \)).
    \end{enumerate}
\end {definition}

\begin{definition}
    Let \( L \) be a positive function on on \( X \). A family \( \mu \) of
    bounded subsets of \( X \) is \emph{controlled by \( L \)} if for any \( A
    \in \mu \) and any \( x \in A \) the diameter of \( A \) does not exceed \(
    L(x) \).
\end{definition}

For a combable geodesic metric space \( X \), dispersed families can be coarsely
homotoped into dispersed sets:

\begin{lemma}
  \label{lem:dispersed-family-set}
  Suppose \( X \) is a bicombable geodesic metric space, and let \( \U \) be a
  dispersed family of subsets of \( X \). Then the union \( U = \bigcup _{a \in \U}
  a \) is coarsely homotopic to a dispersed set in \( X \).
\end{lemma}

\begin{proof}
    For each \( A \in \U \), choose a basepoint \( e(A) \in A \). Since \( X \) is
    bicombable and \( A \) is bounded, we can use lemma \ref{lem:bicomb-property-p}
    to get a homotopy \( H_A : A \times I \to X \) to the point \( e(A) \). We can
    paste these individual maps to get a homotopy \( H: I (U) \to X \) from \( U =
    \bigcup _{A \in \U} A \) to \( V = \bigcup _{A \in \U} \{ e(A) \} \). 
    Since all of these homotopies are \( \rho \)-controlled for the same \( \rho
    \) (property 2 of definition \ref{def:property-p}), our homotopy is coarse.
\end{proof}

\begin{lemma}
    \label{lem:dispersed-categorical}
    Suppose \( X \) is a bicombable proper geodesic metric space which is coarsely
    path-connected, and \( U \subset X \) is a dispersed set. Then \( U \) 
    is coarsely categorical.
\end{lemma}

Before we prove lemma \ref{lem:dispersed-categorical} we need to prove a
technical lemma:

\begin{lemma}
    \label{lem:technical}
    Let \( X \) be a proper geodesic metric space which is coarsely path-connected, 
    and choose a coarse ray \( \alpha : \Rplus \to X \), denote
    \( p = \alpha(0) \). Suppose a subset \( A \subset X \) is dispersed. Then 
    there are functions \( \gamma, \rho : \Rplus \to \Rplus \) having the
    following properties:
    \begin{enumerate}
        \item \( \gamma(R) > R \) for all \( R > 0 \).
        \item Fix any \( R > 0 \). For each point \( x \in A \setminus B _{p} (\gamma(R)) \),
              we can define a path \( h_x : [0, k_x] \to X \) which have the properties:
              \begin{enumerate}
                  \item \( h_x (0) = x, h_x(k_x) \in \im (\alpha) \).
                  \item \( \im (h_x) \cap B _{p} (R) = \emptyset \).
                  \item \( h_x \) is \( \rho \)-controlled.
                  \item \( k_x \leq 2 \norm x \).
              \end{enumerate}
    \end{enumerate}
\end{lemma}

\begin{figure}
    \centering
    \includegraphics{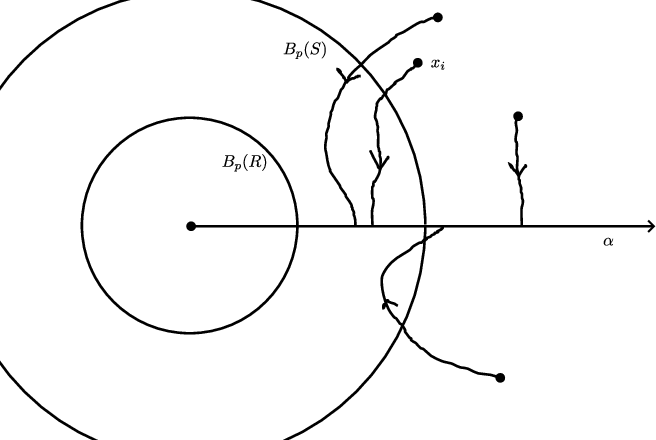}
    \caption{Illustration of lemma \ref{lem:technical}}
    \label{fig:technical-illustration}
\end{figure}

What lemma \ref{lem:technical} states is that for any radius \( R > 0 \), there
is a larger radius \( \gamma(R) = S > 0 \) such that all points in \( A \)
outside the \( S \)-ball around \( p \) can be joined to the coarse ray
\( \alpha \) via a path that does not intersect \( B _{p} (R) \).
Furthermore, we can choose these paths to be nice: uniformly \( \rho
\)-controlled, and at most \( 2 \norm x \)-long. An illustration is given in 
figure \ref{fig:technical-illustration}.

\begin{proof}
  [proof of lemma \ref{lem:technical}]
  First we prove such a path exists, then we prove the stronger conditions of the
  path. Suppose such a path cannot exist. So there exists a radius \( R \) such
  that there exist points of \( A \) arbitrarily far away from \( p \) such that
  all paths joining those points to the ray \( \alpha \) pass through \( B _{p}
  (R) \). Choose a sequence of such points \( \{ x_i \in A \} _{i \in \N} \) such
  that \( r_i = \norm { x_i } \nearrow \infty \).

  Note that for each \( x_i \) for \( i > 1 \), we can choose a geodesic \( g_i
  : [0, r_i] \to X \) such that \( g_i(0) = p \) and \( g_i (r_i) = x_i \). For
  any \( j \in \{ 1, 2, \cdots,  i-1 \} \) define \( x_i ^{j} = g_i (r_j) \), so
  we have \( \norm {x_i^j} = r_j \). For convenience, define \( x_i ^i = x_i \).
  Since \( X \) is a proper metric space, the spheres \( S_p (r_j) \) are
  compact, hence the set \( \{ x_i^j \ | \ i \geq j \} \subset S _p (r_j) \) has
  an accumulation point \( y_j \in S_p (r_j) \). Suppose \( d(y_j, y_{j+1}) =
  d_j \). Construct geodesics \( \beta_j : [0, d_j] \to X \) from \( y_j \) to
  \( y_{j+1} \). Pick a geodesic \( \gamma : [0, r_1] \to X \) from \( p \) to
  \( y_1 \). Now we can define a map \( \beta : \Rplus \to X \) as
  \begin {equation}
    \beta (t) =
    \begin {cases}
        \gamma(t) & \text { if } t \leq r_1 \\
        \beta_j (t - (r_1 + \cdots + r_{j-1})) &
            \text { if } r_1 + \cdots + r_{j-1} < t \leq r_1 + \cdots + r_j
    \end {cases}
  \end {equation}

  Note that \( \beta \) is concatenation of \( \gamma \) and the \( \beta_j
  \)'s. Since \( \beta \) is piecewise geodesic, the map \( \beta \) is
  controlled. It is easy to see that \( \beta \) is proper as well, so \(
  \beta \) is a coarse map.

  Since \( X \) is coarsely path-connected, we can define a coarse
  homotopy \( H : I_d \Rplus \to X \) such that \( H \circ i_0 = \alpha \) and
  \( H \circ i_1 = \beta \). But this is a contradiction!

  Indeed, since \( H \) is a coarse map, in particular it is proper. So there
  exists some \( S > 0 \) such that points \( y \in \im(\beta) \setminus B_{p}
  (S) \) do not touch the ball \( B_{p} (R) \) under the coarse homotopy \( H \).
  We pick a \( r_{j} > S \), and since \( y_{j} \) is an accumulation point in \(
  S_{p} (r_j )\), there exists some \( x_i ^j \in S_{p} (r_j) \) such that \(
  d(y_j, x_i^j) < \varepsilon \) for some \( \varepsilon \ll |R - r_j| \). Let \(
  \varphi : [0, l_j] \to X \) be the path traced by \( y_j \) under the homotopy
  \( H \). Pick a geodesic \( h \) from \( y_j \) to \( x_i^j \). Since \( len(h)
  = \varepsilon \ll |R - r_j| \), this geodesic cannot touch \( B_{p} (R) \).
  Also note that using \( g_i \) we can define a geodesic \( h' \) that goes from
  \( x_i \) to \( x_i^j \) which does not touch \( B_{p} (R) \) either.
  Concatenating these three paths \( h', h, \varphi \) we get a path from \( x_i
  \) to \( \im(\alpha) \) which doesn't intersect \( B_{p} (R) \), which is a
  contradiction.

  Now we show that the path can be chosen nicely.

  Since the homotopy \( H \) is coarse, in particular it is \( \rho
  \)-controlled for some \( \rho : \Rplus \to \Rplus \). Since \( H \) is \(
  \rho \)-controlled, so is \( \varphi \). The maps \( h,h' \) are geodesics,
  so they are \( 1 \)-controlled by definition. Therefore the concatenation is
  \( \rho' \)-controlled for \( \rho' = \max \{ \rho, 1 \} \). This \( \rho' \)
  does not depend on \( i \); all the paths \( h_x \) are controlled by the same
  \( \rho' \).

  For brevity, let \( x = x_i \). \( k_x \) is the length of the concatenation
  \( h_x = \varphi \circ h \circ h' \), therefore \( k_x = len (\varphi) + len
  (h) + len (h') \). Here \( len \) refers to the diameter of the domain of a
  path. So \( len(\varphi) = \norm {y_j} = r_j \). Since \( h \) is a geodesic
  of length \( \varepsilon \), \( len(h) = \varepsilon \), and \( h' \) is a
  geodesic between \( x_i \) and \(x_i^j \), so \( len(h') = r_i - r_j \).
  Adding them,
  \[
       k_x = r_j + \varepsilon + (r_i - r_j) = r_i + \varepsilon
          = \norm x + \varepsilon < 2 \norm x
  \] 
  Since \( \varepsilon \ll |R - r_j| \). This concludes the proof.
\end{proof}

Now we can prove the previous lemma.
\begin{proof}
[proof of lemma \ref{lem:dispersed-categorical}]

Let us fix a basepoint \( p \in X \) and a geodesic ray \( \alpha : \Rplus \to
X\) such that \( \alpha(0) = p \). For convenience, let us define a function \(
\gamma : \Rplus \to \Rplus \), where \( \gamma(0) = 0 \), and for any \( R >
0\), define \( \gamma (R) = R' \) where \( R' \) is as obtained from lemma
\ref{lem:technical}.

Choose an \( R_1 \geq 1 \). Let \( R_2 = \max \{ 2, \gamma(R_1) \}, R_3 = \max
\{ 3, \gamma(R_2) \} \) etc. Using lemma \ref{lem:technical} we will construct
a path \( h_x : [0,k_k] \to X \) for some \( k_x > 0 \) such that
\begin{enumerate}
    \item \( h_x(0) = x, h_x(k_x) \in \im (\alpha) \) for all \( x \in U \).
    \item \( \im (h_x) \cap B _{p} (R_i) = \emptyset\) for all \( x \in U \) such
          that \( \norm x \geq R_{i+1}\).
    \item Each \( h_x \) is \( \rho \)-controlled, for the same \( \rho \).
  \item \( k_x \leq 2 \norm x \).
\end{enumerate}

For \( i=1 \), for each \( x \in U \cap B _{p} (R_1) \), choose any \( h_x \)
satisfying (1). Property 2,3 are satisfied vacuously.

Inductively, suppose we have chosen \( h_x \)'s for all \( x \in U \cap B _{p}
(R_i) \) for \( i = 1, \cdots, k \) which satisfy 1,2,3,4. For each \( x \in U \in
A_p (R_k, R_{k+1})\), by lemma \ref{lem:technical} there exist a path from \(
x\) to \( \im(\alpha) \) which avoids \( B _{p} (R_{k-1}) \) (since \( R_k \geq
\gamma(R_{k-1}) \)). Let \( h_x : [0,1] \to X \) denote that path. By
construction it satisfies all properties 1,2,3 and 4.

Using these paths \( h_x \) for \( x \in A \), we can now construct a coarse
homotopy of \( U \) which is controlled because of 3, and proper because of
2. It will be well-defined because of 4.

% \adi {Should I spell out the construction? It is doable but it will take some
%   space.}

\end{proof}

Now we can combine lemmas \ref{lem:dispersed-family-set} and
\ref{lem:dispersed-categorical} to get the following result:

\begin{proposition}
    Let \( X \) be a bicombable proper geodesic metric space which is coarsely
    path-connected, and let \( \U \) be a dispersed family of subsets of 
    \( X \). Then \( U = \bigcup _{A \in \U} A \) is a coarsely categorical set.
\end{proposition}

\subsection {Asymptotic Dimension} 

Recall that a family of subsets \( \U \) of a metric space \( X \) is called
\emph{uniformly bounded} if there is some \( K > 0 \) for which \( diam(A) \leq
K\) for all \( A \in \U \). The family is \emph { \( r \)-disjoint } if \( d(A,
A') > r \) for every \( A \neq A' \in \U \). Here \( d(A,A') \) is defined to be
\( \inf \{ d(x, x') \ | \ x \in A, x' \in A' \} \). Also we recall the following
definition from \cite{bellAsymptoticDimensionBedlewo2005}:
\begin{definition}
    For a subset \( U \subset X \) and a family of subsets \( \V \) of \( X \),
    the \emph{\( r \)-saturation of \( U \)} is defined as
    \[
        N_r (U, \V) = U \cup \left( \bigcup V \right)
    \]
    Where the big union inside parentheses is over all elements \( V \in \V \)
    such that \( d(U,V) < r \).
\end{definition}

\begin{definition}
    \cite{gromovAsymptoticInvariantsInfinite1993}
  Let \( X \) be a metric space. We say that the \emph{asymptotic dimension} of
  \( X \) does not exceed \( n \) and write \( \asdim X \leq n \) provided for
  every uniformly bounded open cover \( \V \) of \( X \) there is an uniformly
  bounded open cover \( \U \) of \( X \) of multiplicity \( \leq n + 1 \) so
  that \( \V \) refines \( \U \). We write \( \asdim(X) = n \) if it is true
  that \( \asdim (X) \leq n \) and \( \asdim (X) \not \leq n-1 \).
\end{definition}

We are more interested in the following characterization of asymptotic
dimension:

\begin{theorem}
  \label{thm:asdim}
    Let \( X \) be a metric space. The following conditions are equivalent
    \begin{enumerate}
        \item \( \asdim(X) \leq n \);
        \item for every \( r < \infty \) there exist \( r \)-disjoint families \( \U^0,
              \cdots, \U^n \) of uniformly bounded subsets of \( X \) such that \( \bigcup _i
              \U^i \) is a cover of X.
    \end{enumerate}
\end{theorem}

A proof, along with other characterizations can be found in
\cite{bellAsymptoticDimension2008}.

For spaces with finite asymptotic dimension, one can construct dispersed
families of subsets in a nice way, thanks to the following lemma:

\begin{lemma} \label{lem:asdim-dispersed-sets}
    [\cite{DRANISHNIKOV1998791}, lemma 2.9]
    If \( \asdim X \leq n \), then for every positive function \( L \) on \( X
    \) such that \( \lim _{x \to \infty} L(x) = \infty \), there exist
    dispersed families \( \mu_0 , \cdots , \mu_n \) of bounded subsets of \(
    X\) such that each \( \mu_i \) is controlled by \( L \) and \( \mu =
    \bigcup _{i=0} ^n \mu _i \) almost covers \( X \).
\end{lemma}
Here, a family of subsets \emph{almost covers} \( X \setminus K \) for some
compact \( K \subset X \).

Now we have enough machinery to prove our main result.
\begin{theorem}
  \label{thm:ccat-asdim}
  Let \( X \) be a bicombable proper geodesic metric space which is coarsely
    path-connected. Then we have the following inequality 
    \[
      \ccat (X) \leq \asdim(X).
  \]
\end{theorem}

\begin{proof}
  Let us fix a basepoint \( p \in X \). Let us assume that \( \asdim (X) \leq n
  \). In light of lemma \ref{lem:dispersed-categorical} it is enough to
  construct \( n+1 \) families of subsets \( \{ \U^i \}_{i=0} ^{i=n} \) such
  that
  \begin{enumerate}
    \item Each \( \U^i \) is a dispersed family of subsets, and
    \item \( X \) is covered by these \( \U^i \)'s:
      \[
        \bigcup _{i=0} ^{n} \left(
          \bigcup _{A \in \U^i} A
        \right) = X.
      \]
  \end{enumerate}

  One way of constructing these \( \U^i \)'s is using lemma
  \ref{lem:asdim-dispersed-sets}, we can take our positive function \( L \) to be the
  distance function.

  We describe an alternative way of constructing these sets in the appendix.
\end{proof}

The above theorem can be stated for bicombable groups now:

\begin{theorem} \label{thm:ccat-asdim-grp}
    For a bicombable group \( \Gamma \) that is coarsely path-connected, we have the
    inequality
    \[
        \ccat(\Gamma) \leq \asdim (\Gamma)
    .\]
\end{theorem}

And as a result of claim \ref{clm:groups-semistable-pi0}, we have the following

\begin{corollary} \label{cor:ccat-asdim-grp-ends}
    For any bicombable 1-ended discrete group \( \Gamma \) which is semistable at
    \( \infty \) we have the inequality
    \[ 
        \ccat(\Gamma) \leq \asdim (\Gamma) 
    .\]
\end{corollary}

\section {Acknowledgements}
The author would like to thank his advisor Alexander Dranishnikov for all his
help and insightful discussions, and Thomas Weighill for helping him understand
proper vs coarse homotopy. The author would also like to thank UF Center for
Applied Mathematics for summer fellowship award.

\printbibliography

\appendix
\newpage
\section {Appendix: Construction of dispersed sets}

Given a metric space \( X \) with \( \asdim X \leq n \), we describe how to
construct \( n+1 \) families of subsets \( \{ \U^i \}_{i=0} ^{i=n} \) such
that
\begin{enumerate}
\item Each \( \U^i \) is a dispersed family of subsets, and
\item \( X \) is covered by these \( \U^i \)'s:
  \[
    \bigcup _{i=0} ^{n} \left(
      \bigcup _{A \in \U^i} A
    \right) = X.
  \]
\end{enumerate}

We inductively define an increasing sequence of radii \( \{ R_j \}_{j \in \N}
\), scales \( \{ \lambda_j \} _{j \in \N} \), and families of subsets \( \C^i_j
\) inside \( A_p (R_{j-1}, R_j) \) for \( j > 0 \). For convenience, let \( R_0
= \lambda_0 = 0 \).

For the base case \( j = 1 \), choose a scale \( \lambda_1 > 0 \). Since \(
\asdim (X) \leq n \), using \ref{thm:asdim} we can get uniformly \( D_1
\)-bounded, \( \lambda_1 \)-disjoint families of subsets \( \B_1 ^{0}, \B_1
^{1}, \cdots , \B_1 ^{n} \) for some \( D_1 > 0 \) which cover \( X \). Choose
\( R_1 \gg \max \{ \lambda_1, D_1 \} \). Now for \( i = 0, 1, \cdots, n \)
define 
\[
  \C_1 ^{i} = \{ B \cap B _{p} (R_1) \ | \ B \in \B_1 ^{i} \}.
\]
This concludes our base step. Now inductively assume that we have already
constructed \( R_j, \C_j ^{i}, \lambda_j \) for \( j = 1, 2, \cdots, k \).
Choose a scale \( \lambda _{k+1} \gg R_k \). Using \ref{thm:asdim}, we get
uniformly \( D_{k+1} \)-bounded \( \lambda_{k+1} \)-disjoint families of
subsets \( \B _{k+1} ^0, \B_{k+1}^1, \cdots, \B_{k+1}^n \) which cover \( X
\), for some \( D_{k+1} > 0 \). We define
\[
C_{k+1} ^{i} = \{ B \cap A_p (R_k, R_{k+1}) \ | \ B \in \B _{k+1}^{i} \}.
\]

\begin{figure}
    \centering
    \includegraphics{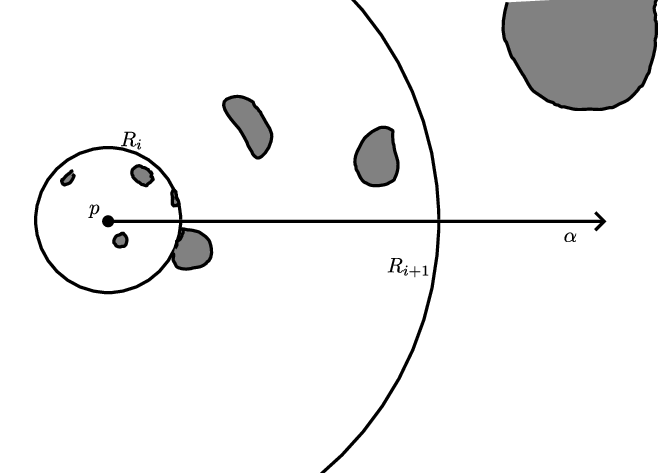}
    \caption{Illustration of \( C_{k+1}^{i} \). The shaded blobs are the
    elements of \( C_{k+1}^{i} \).}
    \label {fig:3}
\end{figure}

Figure \ref{fig:3} shows this construction.

By our construction, \( \bigcup_j \C_j ^i \) is not quite a dispersed family,
as sets in \( \C_j^{i} \) can be very close to some set in \( \C_{j+1}^i \).
We need to redistribute the sets slightly to get rid of these possibilities.

For \( j=1 \) define
\[
   \D_1^{i} = \{ C \in \C_1^i \ | \ d(C, C') \geq \lambda_1 
   \text { for all } C' \in \C _2 ^{i} \}.
\]

For \( j > 1 \) define
\[
D_j ^{i} = \{ N _{\lambda_{j-1}} (C, \C_{j-1} ^{i}) \ | \ 
  C \in \C_j ^{i}, \ d(C, C') \geq \lambda_j 
  \text { for all } C \in \C_{j+1} ^{i} \}.
\]

Finally, we can define \( \U^{i} = \D_1^{i} \cup \D_2^{i} \cup \D_3 ^{i} \cup
\cdots \).
\begin{itemize}

\item Each family \( \D^i = \{ \D^i_1, \D^i_2, \cdots \} \) consists of
        pairwise disjoint sets. This follows from the fact that \( \lambda_{j+1} \gg
  \lambda _j\), so an element of \( \C_j ^{i} \) cannot be \( \lambda_j
  \)-close to two or more elements in \( \C_{j+1}^{i} \). Boundedness
  follows from the fact that \( R_j \gg \max \{ D_j, \lambda_j \} \), so
  elements in \( \C_j ^{i} \) cannot be simultaneously \( \lambda_{j-1}
  \)-close to an element in \( \C_{j-1}^{i} \) and \( \lambda_j \)-close to
  an element in \( \C_{j+1} \).

    \item Each family \( \D^{i} \) is a dispersed family of subsets. This
        follows from the fact that each set in \( \D^i_j \) is \( \lambda_j
        \)-disjoint from others, and from sets in \( \D^{i}_{j+1} \). Since
        element of \( \D^{i} \) is bounded, we can conclude that the function 
        \[
            \partial (r) = \inf \{ d(x_1, x_2) 
                    : x_1 \neq x_2 \in U \setminus B _{p} (r) \}
        \]
        goes to infinity.
\end{itemize}
This concludes our construction.
\qed
\end {document}